# A FUNCTIONAL CLT FOR THE OCCUPATION TIME OF A STATE-DEPENDENT BRANCHING RANDOM WALK

By Matthias Birkner and Iljana Zähle[1]

*Weierstrass-Institut für Angewandte Analysis and Stochastik
and Universität Erlangen–Nürnberg*

We show that the centred occupation time process of the origin of a system of critical binary branching random walks in dimension $d \geq 3$, started off either from a Poisson field or in equilibrium, when suitably normalized, converges to a Brownian motion in $d \geq 4$. In $d = 3$, the limit process is a fractional Brownian motion with Hurst parameter $3/4$ when starting in equilibrium, and a related Gaussian process when starting from a Poisson field. For (dependent) branching random walks with state dependent branching rate we obtain convergence in f.d.d. to the same limit process, and for $d = 3$ also a functional limit theorem.

**1. Introduction and main result.** We study the fluctuation behavior of the occupation time in a single point of a system of critical binary branching random walks with a state-dependent branching rate (BRW). BRW consists of particles which move independently on $\mathbb{Z}^d$ in continuous time according to a given random walk kernel $a$. The branching rate at a site depends on the number of particles there: if there are $k$ individuals at $x$, a branching event in $x$ occurs at rate $\sigma(k)$. The particle which is chosen to branch then leaves either two or zero offspring at its current location, each possibility occurring with probability $1/2$. The classical case of independent branching with constant branching rate $\rho$ corresponds to $\sigma(k) = \rho k$. We assume further on that the branching rate function $\sigma$ is not $\equiv 0$ and Lipschitz:

$$(1.1) \qquad |\sigma(m) - \sigma(n)| \leq c|m - n|,$$

Received April 2005; revised September 2006.
[1]Supported in part by DFG Schwerpunkt "Interagierende stochastische Systeme von hoher Komplexität."

which especially implies

(1.2) $$\sigma(k) \leq c_2 k \quad \text{for some } 0 < c_2 < \infty.$$

We denote by $\xi_t(x)$ the number of particles present at location $x$ at time $t$. We assume that the transition rate matrix $a(x,y) = a(0, y-x)$ governing the individual motion of particles is symmetric, irreducible and has finite second moments, which implies

(1.3) $$(Q_{ij})_{i,j=1,\ldots,d} = \left( \sum_{x \in \mathbb{Z}^d} a(0,x) x_i x_j \right)_{i,j=1,\ldots,d} \quad \text{is finite and invertible},$$

where $x = (x_i)_{i=1,\ldots,d}$. We have $\sum_x a(0,x)x = 0$ by symmetry, and we can assume without loss of generality that $a$ is stochastic, that is, $\sum_x a(0,x) = 1$. Denote the continuous time transition probabilities by $a_t(x,y)$.

It is well known that the independent BRW in $d \leq 2$, starting from any initial condition with bounded local density, suffers local extinction, that is, $\xi_t(x) \to 0$ in probability as $t \to \infty$ for any $x \in \mathbb{Z}^d$. For translation invariant, ergodic initial distributions with finite intensity, this can be found in, for example, [17]. It can be extended by the comparison argument in the proof of (5.1) in [13].

On the other hand, in $d \geq 3$, there exists a one-parameter family of extremal invariant probability measures $\Lambda_\vartheta$, $\vartheta \geq 0$, parametrized by the expected density: $\int \xi(x) \Lambda_\vartheta(d\xi) = \vartheta$. Each $\Lambda_\vartheta$ is shift-invariant, and $\{\xi_t(x) : x \in \mathbb{Z}^d, t \geq 0\}$ under $\Lambda_\vartheta$ is ergodic with respect to space- and time-shifts. See references below Theorem 2.3.

Let us denote the distribution of a Poisson field on $\mathbb{Z}^d$ with homogeneous intensity $\vartheta$ by $\mathcal{H}(\vartheta)$, that is, under $\mathcal{H}(\vartheta)$, the random variables $\xi(x)$, $x \in \mathbb{Z}^d$, are i.i.d. Poisson$(\vartheta)$. If $\mathcal{L}(\xi_0) = \mathcal{H}(\vartheta)$, we have $\mathcal{L}(\xi_t) \to \Lambda_\vartheta$ weakly as $t \to \infty$.

First we present the main result for the case of independent branching [i.e., $\sigma(k) = \rho k$]. Let $\mathcal{L}(\xi_0) \in \{\mathcal{H}(\vartheta), \Lambda_\vartheta\}$. By ergodicity, the occupation time of any point $x \in \mathbb{Z}^d$ satisfies

$$\frac{1}{T} \int_0^T \xi_t(x) \, dt \to \vartheta \quad \text{almost surely as } T \to \infty.$$

Thus, a natural question concerns the random fluctuations of the occupation time around its asymptotic limit. This is the content of our main result:

THEOREM 1.1 (Independent branching). 1. *If $(\xi_s)_{s \geq 0}$ is started in the (unique extremal) equilibrium distribution $\Lambda_\vartheta$ with intensity $\vartheta > 0$, then the processes*

(1.4) $$X_t^N := \frac{1}{h_d(N)} \int_0^{Nt} (\xi_s(0) - \vartheta) \, ds, \quad t \geq 0,$$

*converge in distribution as elements of $\mathcal{C}([0,\infty),\mathbb{R})$ (denoting the set of continuous functions on $[0,\infty)$ with values in $\mathbb{R}$) toward a Brownian motion in $d \geq 4$ and to a fractional Brownian motion with Hurst parameter $3/4$ in $d = 3$ as $N \to \infty$, where the norming is given by*

$$(1.5) \qquad h_d(t) = \begin{cases} t^{3/4}, & d=3, \\ \sqrt{t \log t}, & d=4, \\ \sqrt{t}, & d \geq 5. \end{cases}$$

*The covariance of the limiting process $X$ is given by*

$$\operatorname{Cov}(X_s, X_t)$$
$$= \begin{cases} \dfrac{\sqrt{2}}{3\pi^{3/2}} (\det Q)^{-1/2} \vartheta \rho [t^{3/2} + s^{3/2} - |t-s|^{3/2}], & d=3, \\ (2\pi)^{-2} (\det Q)^{-1/2} \vartheta \rho \times (s \wedge t), & d=4, \\ \left[ 2 \int_0^\infty du\, a_u(0,0) + \rho \int_0^\infty du\, u a_u(0,0) \right] \vartheta \times (s \wedge t), & d \geq 5. \end{cases}$$

2. *The same conclusions hold if $\mathcal{L}(\xi_0) = \mathcal{H}(\vartheta)$, and $d \geq 4$. In the case $\mathcal{L}(\xi_0) = \mathcal{H}(\vartheta)$ and $d = 3$, the processes $X^N$ converge toward a Gaussian process $X$ with covariance given by*

$$(1.6) \qquad \operatorname{Cov}(X_s, X_t) = \dfrac{2\sqrt{2}}{3\pi^{3/2}} (\det Q)^{-1/2} \vartheta \rho$$
$$\times \left[ t^{3/2} + s^{3/2} - \dfrac{1}{2}|t-s|^{3/2} - \dfrac{1}{2}(t+s)^{3/2} \right].$$

The normalizations $h_d$ are dictated by the requirement of a nontrivial covariance function for the limit process, and this in turn is determined by the decay properties of the transition probabilities of the underlying random walk $a$; see the calculations in Section 3.3. Note that with $\rho = 0$, BRW becomes a system of independent random walks, and has the family $\mathcal{H}(\vartheta)$, $\vartheta \geq 0$, of shift-invariant extremal equilibria. In the situation $\rho = 0$, we see from Theorem 1.1 that under the normalizations used in Theorem 1.1 the limit process $X$ is trivial in $d \leq 4$ and a Brownian motion in $d \geq 5$. This is in keeping with a "metatheorem" that the introduction of branching shifts "critical dimensions" by 2: In a system of independent random walks, the occupation time requires normalization by $t^{3/4}$ in $d = 1$, $\sqrt{t \log t}$ in $d = 2$ and $\sqrt{t}$ in $d \geq 3$ in order to obtain a nontrivial limit (see [12]).

While for nonbranching random walks, the nonclassical norming is due to recurrence properties of the individual particles, the behavior in our case is governed by the recurrence properties of *families*: The equilibrium of a BRW can be decomposed into a Poisson system of "clans" of particles with a common ancestor (see, e.g., [28]), and such a clan will visit the origin infinitely often if and only if $d \leq 4$ (see [26], Theorem 1). This allows in the case of independent branching, at least on a heuristic level, also to understand the different normings. Substituting the probabilistic representation of the Palm distribution at 0 at time $T$ from [26], Proposition 1 for "a typical clan which contributes to the occupation time up to $T$," we see that the expected contribution per clan is

$$\text{(1.7)} \qquad \int_0^T dt \int_A^t ds\, a_{T+t-2s}(0,0),$$

where $A = 0$ when $\mathcal{L}(\xi_0) = \mathcal{H}(\vartheta)$ and $A = -\infty$ when $\mathcal{L}(\xi_0) = \Lambda_\vartheta$. In $d=3$ this grows like $\sqrt{T}$, hence, we expect of the order $T/\sqrt{T} = \sqrt{T}$ clans to contribute. So due to independence of clans, the fluctuations should indeed be of order $\sqrt{T}(\sqrt{T})^{1/2} = T^{3/4}$. In $d=4$, this grows like $\log T$, so $T/\log T$ clans should contribute, suggesting fluctuations of the order $\log T \sqrt{T/\log T} = \sqrt{T \log T}$. In $d \geq 5$, (1.7) is bounded, corroborating the classical norming.

It is remarkable that the correlations introduced by the branching are strong enough that in $d=3$, the limit process itself depends on the initial condition, not only on its density. Even though $\xi_t$, starting from $\mathcal{H}(\vartheta)$, converges in distribution to $\Lambda_\vartheta$, the "building up" of equilibrium is reflected in the different covariance structure of the renormalized occupation time process.

Note that the centred Gaussian process $(X_t)$ with covariance given by (1.6) can be represented as $X_t = (B_t^{(3/4)} + B_{-t}^{(3/4)})/\sqrt{2}$, where $(B_t^{(3/4)})_{t \in \mathbb{R}}$ is a fractional Brownian motion with Hurst parameter $3/4$ and $B_0^{(3/4)} = 0$ (see [3]). On the level of variances, this can be seen as follows. By (3.14) and (3.18),

$$\text{Var}^{\Lambda_\vartheta}(X_t^N) - \text{Var}^{\mathcal{H}(\vartheta)}(X_t^N) = \frac{\rho}{2N^{3/2}} \int_0^{Nt} du \int_0^{Nt} dv \int_{u+v}^\infty dr\, a_r(0,0).$$

Lemma 3.3 implies that the right-hand side is equal to $\text{Cov}^{\Lambda(\vartheta)}(X_{-t}^N, X_t^N)$ up to negligible terms. Hence,

$$\text{Var}^{\mathcal{H}(\vartheta)}(X_t^N) \approx \text{Var}^{\Lambda(\vartheta)}\left(\frac{1}{\sqrt{2}}(X_t^N + X_{-t}^N)\right).$$

It remains an intriguing question to explain this representation from the point of view of building up a family structure in the branching random walk.

For state dependent branching in $d \geq 3$, Proposition 3 in [8] shows that

$$\text{(1.8)} \qquad \Lambda_\vartheta = \underset{t \to \infty}{\text{w-lim}} \mathcal{L}^{\mathcal{H}(\vartheta)}(\xi_t)$$

exists for any $\vartheta \geq 0$. $\Lambda_\vartheta$ is a shift-invariant equilibrium and satisfies $\mathbb{E}^{\Lambda_\vartheta}[\xi_0(x)] = \vartheta$ and $\mathbb{E}^{\Lambda_\vartheta}[(\xi_0(x))^2] < \infty$. We denote $\sigma_\vartheta^{\text{eq}} := \mathbb{E}^{\Lambda_\vartheta}[\sigma(\xi_0(x))]$, which is independent of $x$ because of the shift-invariance of $\Lambda_\vartheta$. By (1.2) and the assumption $\sigma \not\equiv 0$, we have $\sigma_\vartheta^{\text{eq}} \in (0, \infty)$.

In the case $d = 3$ with start in the Poisson field we can prove a functional CLT for the occupation time in zero. For $d \geq 4$, we show f.d.d.-convergence of the renormalized occupation time.

THEOREM 1.2 (State dependent branching). 1. *In $d = 3$ the processes $(X_t^N)_{t \geq 0}$ defined in (1.4) with $\mathcal{L}(\xi_0) = \mathcal{H}(\vartheta)$ and with the norming given in (1.5) converge in distribution as elements of $\mathcal{C}([0, \infty), \mathbb{R})$ toward a Gaussian process $X$ with covariance given by*

$$\text{Cov}(X_s, X_t) = \frac{2\sqrt{2}}{3\pi^{3/2}} (\det Q)^{-1/2} \sigma_\vartheta^{\text{eq}}$$
$$\times \left[ t^{3/2} + s^{3/2} - \frac{1}{2}|t - s|^{3/2} - \frac{1}{2}(t + s)^{3/2} \right]. \quad (1.9)$$

2. *In $d \geq 4$ the processes $(X_t^N)_{t \geq 0}$ defined in (1.4) with $\mathcal{L}(\xi_0) \in \{\mathcal{H}(\vartheta), \Lambda_\vartheta\}$ and with the norming given in (1.5) converge in f.d.d.-sense toward a Brownian motion $X$ with covariance*

$$\text{Cov}(X_s, X_t)$$
$$= \begin{cases} (2\pi)^{-2} (\det Q)^{-1/2} \sigma_\vartheta^{\text{eq}} \times (s \wedge t), & d = 4, \\ \left[ 2\vartheta \int_0^\infty du\, a_u(0, 0) + \sigma_\vartheta^{\text{eq}} \int_0^\infty du\, u a_u(0, 0) \right] \times (s \wedge t), & d \geq 5. \end{cases} \quad (1.10)$$

REMARK 1.1. For the state dependent branching case, we can at present not prove tightness in $d \geq 4$ due to a lack of manageable higher moment formulas. We conjecture that in the case $d = 3$, when starting in equilibrium, the limit process would again be fractional Brownian motion. In order to prove this using the techniques employed in the present paper, we would require an equivalent of the main result from [28], namely, a spatial renormalization result for the equilibrium, in the state dependent case (see the proof of Lemma 3.2). While we believe this to be true, the techniques used in [28] depend on infinite divisibility, and can thus not be readily adapted.

Corresponding functional central limit theorems for the occupation time of *reversible* interacting particle systems are well known; see, for example, [19, 24], or, more generally, [20] for central limit theorems for additive functionals of reversible Markov processes. In the nonreversible situation of independently branching system, nonfunctional versions of central limit theorems have been obtained in [14].

In order to obtain tightness in $d \geq 4$ (independent case), we have traded reversibility for infinite divisibility, which opens the possibility of rather explicit calculations. 4th moment calculations are feasible, although cumbersome, because of the independence of families founded by different particles.

A program similar to ours has been carried out by Bojdecki, Gorostiza and Talarczyk in [6] and [7] in the case of independent branching in a somewhat different scenario with completely different techniques: They consider critical binary branching particles in $\mathbb{R}^d$, where the individual particle moves according to a symmetric $\alpha$-stable process, with $\alpha \in (0, 2]$, and obtain the following results: for $\alpha < d < 2\alpha$, starting from a homogeneous Poisson process, the occupation time requires a nonclassical norming and converges to sub-fractional Brownian motion [the centered Gaussian process with covariance given by (1.6)], whereas the limit process is Brownian for $d \geq 2\alpha$ (with a logarithmic correction to the norming in the boundary case $d = 2\alpha$). Bojdecki, Gorostiza and Talarczyk have also considered the scenario with a heavy-tailed offspring distribution, see [4] and [5].

Our set-up is different in the following respect: we allow for a state dependent branching rate, we consider the lattice instead of continuous space, and we focus on the occupation time of a single point, whereas Bojdecki, Gorostiza and Talarczyk consider the integral of the particle system against test functions from $\mathcal{S}(\mathbb{R}^d)$. As to the techniques, Bojdecki, Gorostiza and Talarczyk rely on computations of Laplace functionals and Fourier analysis, while, in our case, the discreteness of space allows to use martingale decompositions of the occupation time, and to employ techniques from the field of interacting particle systems (similar to [19] and [24]). Our scenario, namely, individual motion with a finite second moment, combined with critical binary branching, corresponds to the case $\alpha = 2$. This invites to conjecture that if we used an individual motion which is in the domain of attraction of an $\alpha$-stable law (with general $\alpha \in (0, 2]$), we would find the same $\alpha$-dependence of regimes as Bojdecki, Gorostiza and Talarczyk. On the other hand, our Theorem 1.1, part 1 suggests that in the scenario of [6], starting off from an extremal equilibrium for the branching system instead of a Poisson process, the limit process should be a fractional Brownian motion. This has, in fact, been proved by Miłoś, [23], using techniques similar to [6].

The rest of this paper is organized as follows: We collect some well-known facts about random walks and branching random walks in Section 2. Convergence and asymptotic Gaussianity of finite-dimensional distributions is proved in Section 3: in the case $d \geq 4$ we decompose the occupation time into a martingale plus an asymptotically negligible remainder term (Section 3.1), in the case $d = 3$ we "distill" a white noise out of the fluctuations of the particle system and represent the occupation time as an integral with respect to this noise (Section 3.2). In order to prove tightness, we use moment estimates.

## 2. Preliminaries.

2.1. *Formulas related to random walks.* The underlying motion process has generator

$$Lf(x) = \sum_{y \in \mathbb{Z}^d} a(x,y)(f(y) - f(x)).$$

The continuous time transition probabilities $a_t(x,y)$ solve the backward equation $\frac{\partial}{\partial t} a_t(x,y) = L a_t(\cdot, y)(x)$. We denote the transition semigroup by $S_t f(x) := \sum_y a_t(x,y) f(y)$. Let

$$g(x,y) := \int_0^\infty a_t(x,y)\,dt$$

be the Green function and

$$g_\lambda(x,y) := \int_0^\infty e^{-\lambda t} a_t(x,y)\,dt$$

the resolvent. We denote the Green operator by $\mathcal{G} f(x) := \sum_y g(x,y) f(y)$. The function $x \mapsto g(x,0)$ is a solution of $-L\phi = \delta_0$ and $x \mapsto g_\lambda(x,0)$ a solution of $\lambda \phi - L\phi = \delta_0$. Define

$$u_t(x,y) := \int_0^t a_s(x,y)\,ds,$$

the Green function of the random walk killed at time $t$. The function $(t,x) \mapsto u_t(x,0)$ solves $(\partial_t - L)\phi = \delta_0, \phi_0(x) \equiv 0$. The Dirichlet form of the underlying random walk is $\sum_{x,y} a(x,y)(\phi(y) - \phi(x))^2 = 2\langle \phi, (-L)\phi \rangle$ for $\phi \in \ell_2(\mathbb{Z}^d)$.

Note that our assumptions on $a$ imply the local CLT

$$(2.1) \quad \sup_{x \in \mathbb{Z}^d} \left( \left( \frac{|x|}{\sqrt{t}} \right)^s + 1 \right) \left| a_t(0,x) - p_t(0,x) \left[ 1 + \sum_{k=1}^{s-2} t^{-k/2} P_k\left( \frac{x}{\sqrt{t}} \right) \right] \right|$$
$$= o(t^{-(d+s-2)/2}),$$

where $P_k$ is a polynomial of degree $3k$ and

$$p_t(0,x) = (2\pi t)^{-d/2} (\det Q)^{-1/2} \exp\left( -\frac{x^T Q^{-1} x}{2t} \right).$$

The local CLT for discrete time random walks can be found in [10] as Corollary 22.3. From that one can derive the corresponding result for continuous time. This can be done similarly to [1], page 113, where a result on the Galton–Watson process is transferred from discrete time to continuous time. Specifically, we need the following form of the local CLT:

$$(2.2) \quad a_t(0,0) = (2\pi t)^{-d/2} \det(Q)^{-1/2} + o(t^{-d/2}) \quad \text{as } t \to \infty.$$

From this, one can conclude that $\|g(\cdot,0)\|_2^2 = \int_0^\infty ds \int_0^\infty dt\, a_{s+t}(0,0)$, so that $g(\cdot,0) \in \ell_2(\mathbb{Z}^d)$ in case $d \geq 5$, whereas $\|g_{1/N}(\cdot,0)\|_2^2 \sim C \log N$ in case $d = 4$.

2.2. *Basic results on branching random walk.* A convenient choice of the state space for a branching random walk (as well as many other "spatially homogeneous" particle systems), going back to Liggett and Spitzer [22], is

$$\mathfrak{X} = \left\{\mu \text{ an integer-valued measure on } \mathbb{Z}^d \colon \sum_{x \in \mathbb{Z}^d} \gamma(x)\mu(x) < \infty\right\},$$

where $\gamma$ is a strictly positive function on $\mathbb{Z}^d$ satisfying $\sum_{y \in \mathbb{Z}^d} a(x,y)\gamma(y) \leq M\gamma(x)$ for some constant $M > 0$. Note that the dependence of $\mathfrak{X}$ on the particular choice of $\gamma$ is irrelevant for our purposes, as any random $(\xi(x))_{x \in \mathbb{Z}^d}$ satisfying $\sup_x \mathbb{E}[\xi(x)] < \infty$ automatically has $\mathbb{P}(\xi \in \mathfrak{X}) = 1$ irrespective of $\gamma$. A formal construction of the independent BRW $(\xi_t)_{t \geq 0}$ as an $\mathfrak{X}$-valued Markov process can be found, for example, in Section 1 of [17]. The BRW with state dependent branching rate is constructed in Section 2.2 of [8]. The generator is given by

$$\mathcal{L}F(\xi) = \sum_{x \in \mathbb{Z}^d} \sum_{y \in \mathbb{Z}^d} \xi(x)a(x,y)(F(\xi^{x,y}) - F(\xi))$$

(2.3)
$$+ \sum_{x \in \mathbb{Z}^d} \frac{\sigma(\xi(x))}{2}(F(\xi^{x,+}) + F(\xi^{x,-}) - 2F(\xi))$$

with $\xi^{x,y} = \xi - \delta_x + \delta_y$, $\xi^{x,+} = \xi + \delta_x$ and $\xi^{x,-} = \xi - \delta_x$. The branching random walk with state dependent branching rate $(\xi_t)_{t \geq 0}$ with initial condition $\xi_0 \in \mathfrak{X}$ can be constructed as the unique solution to

$$\xi_t(x) = \xi_0(x) + \sum_{y \neq x} \left[\int_{[0,t] \times \mathbb{N}} \mathbf{1}(\xi_{s-}(y) \geq n)\bar{N}^{y,x}(ds\,dn)\right.$$

$$\left. - \int_{[0,t] \times \mathbb{N}} \mathbf{1}(\xi_{s-}(x) \geq n)\bar{N}^{x,y}(ds\,dn)\right]$$

$$+ \int_{[0,t] \times \mathbb{N} \times [0,1]} \mathbf{1}\left(\xi_{s-}(x) \geq n, \frac{\sigma(\xi_{s-}(x))}{c_2\xi_{s-}(x)} \geq u\right)\bar{N}^{x,+}(ds\,dn\,du)$$

$$- \int_{[0,t] \times \mathbb{N} \times [0,1]} \mathbf{1}\left(\xi_{s-}(x) \geq n, \frac{\sigma(\xi_{s-}(x))}{c_2\xi_{s-}(x)} \geq u\right)\bar{N}^{x,-}(ds\,dn\,du)$$

for all $x \in \mathbb{Z}^d$, $t \geq 0$. Here, $\bar{N}^{x,y}$, $x \neq y$, are independent Poisson processes on $[0,\infty) \times \mathbb{N}$ and $\bar{N}^{x,+}, \bar{N}^{x,-}$, $x \in \mathbb{Z}^d$, are independent Poisson processes on $[0,\infty) \times \mathbb{N} \times [0,1]$, all independent of $\xi_0$. $\bar{N}^{x,y}$ has intensity measure $a(x,y)\,dt \otimes d\ell$, $\bar{N}^{x,+}, \bar{N}^{x,-}$ have intensity measure $(c_2/2)\,dt \otimes d\ell \otimes du$ ($dt$, $du$ are Lebesgue measures, $\ell$ is counting measure). For fixed $\xi_0$, $(\xi_t)$ is adapted to the filtration generated by these Poisson processes. See, for example, [8],

Lemma 1 and Remark 3. Define
$$N_t^{x,y} := \int_{[0,t]\times\mathbb{N}} \mathbf{1}(\xi_{s-}(x) \geq n)\bar{N}^{x,y}(ds\,dn),$$
$$N_t^{x,\pm} := \int_{[0,t]\times\mathbb{N}\times[0,1]} \mathbf{1}\left(\xi_{s-}(x) \geq n, \frac{\sigma(\xi_{s-}(x))}{c_2\xi_{s-}(x)} \geq u\right)\bar{N}^{x,\pm}(ds\,dn\,du)$$

(with the obvious interpretations: $N^{x,+}$ counts the number of births at $x$, $N^{x,-}$ counts the number of deaths at $x$, $N^{x,y}$ counts how many times a particle jumps from $x$ to $y$). Thus, we can rewrite
$$\xi_t(x) = \xi_0(x) + N_t^{x,+} - N_t^{x,-} + \sum_{y \neq x}(N_t^{y,x} - N_t^{x,y}), \qquad x \in \mathbb{Z}^d, t \geq 0.$$

Immediately from the independence properties of the driving Poisson processes $\bar{N}$ we get the following:

LEMMA 2.1. *Assume that* $\sup_x \mathbb{E}[\xi_0(x)^2] < \infty$. *The compensated processes*
$$\tilde{N}_t^{x,y} := N_t^{x,y} - a(x,y)\int_0^t \xi_s(x)\,ds,$$
$$\tilde{N}_t^{x,\pm} := N_t^{x,\pm} - \tfrac{1}{2}\int_0^t \sigma(\xi_s(x))\,ds$$
*are pairwise orthogonal, square integrable martingales with compensators given by*
$$\langle\tilde{N}^{x,y}\rangle_t = a(x,y)\int_0^t \xi_s(x)\,ds, \langle\tilde{N}^{x,+}\rangle_t = \langle\tilde{N}^{x,-}\rangle_t = \tfrac{1}{2}\int_0^t \sigma(\xi_s(x))\,ds.$$

For $f_t \in \ell_1(\mathbb{Z}^d)$, put
$$(2.4) \qquad F_t(\xi) := \langle f_t, \xi - \vartheta\lambda\rangle = \sum_{x\in\mathbb{Z}^d} f_t(x)(\xi(x) - \vartheta).$$

Note that this sum is well defined if $\sup_x \mathbb{E}|\xi(x) - \vartheta| < \infty$.

By compensating the driving Poisson processes, we obtain the following:

LEMMA 2.2. *Let* $f:[0,\infty)\times\mathbb{Z}^d \to \mathbb{R}$ *satisfy* $\sup_{t\leq T}(\|f_t\|_1 + \|\partial_t f_t\|_1 + \|Lf_t\|_1) < \infty$, *and let* $F_t$ *be defined by (2.4). Then we have*
$$(\partial_t + \mathcal{L})F_t(\xi) = \langle(\partial_t + L)f_t, \xi - \vartheta\ell\rangle,$$
*for* $t \in [0,T]$, *and the martingale* $M_t := F_t(\xi_t) - F_0(\xi_0) - \int_0^t(\partial_s + \mathcal{L})F_s(\xi_s)\,ds$, $0 \leq t \leq T$, *can be represented as*
$$M_t = \sum_{x\in\mathbb{Z}^d}\int_0^t f_s(x)(d\tilde{N}_s^{x,+} - d\tilde{N}_s^{x,-}) + \sum_{x,y\in\mathbb{Z}^d}\int_0^t(f_s(y) - f_s(x))\,d\tilde{N}_s^{x,y}.$$

We state some basic properties of a critical (finite variance) branching random walk in $d \geq 3$ which we will need in the following. Let $\hat{a}(x,y) = \frac{1}{2}(a(x,y) + a(y,x))$ be the symmetrized transition kernel. (In our case $\hat{a} = a$.)

PROPOSITION 2.3. *Assume that $\hat{a}$ is transient. Then for each $\vartheta \geq 0$, there exists an extremal invariant probability measure $\Lambda_\vartheta \in \mathcal{P}(\mathcal{N}(\mathbb{Z}^d))$ with $\int \xi(0) \Lambda_\vartheta(d\xi) = \vartheta$. Each $\Lambda_\vartheta$ is translation invariant. If $\mathcal{L}(\xi_0) \in \{\mathcal{H}(\vartheta), \Lambda_\vartheta\}$,*

$$\frac{1}{t} \int_0^t f(\xi_s) \, ds \xrightarrow[t\to\infty]{} \int f(\xi) \Lambda_\vartheta(d\xi) \qquad in \ L_1$$

*for linear bounded, local functions $f$.*

REMARK 2.1. In case of independent branching $\Lambda_\vartheta$ is unique and the convergence also holds almost surely.

For the independent branching case, the earliest version, in a discrete generations setting, appeared in [21], Satz 5.4 and Satz 5.5. A corresponding statement for a continuous-time model (which differs from our definition of a branching random walk only insofar as birth and motion are coupled) is contained in Theorem 6.3 and Theorem 6.4 of [15]. A proof of Proposition 2.3 can also be obtained from the proof of Theorem 2(a), Case 1 in [17] by specializing to $p = 0$. Results in this spirit are well known, see, for example, the references given for Theorem 0 in [2], which states the corresponding results for two "continuous relatives" of branching random walk, namely, branching Brownian motion and the Dawson–Watanabe superprocess.

For Proposition 2.3 in the state dependent branching case, we refer to Proposition 3 in [8], where $\Lambda_\vartheta$ is constructed as the weak limit of the system started in the Poisson field, and Theorem 1 in [9], which proves convergence of time averages.

### 3. Finite-dimensional distributions.

PROPOSITION 3.1. *Let $\mathcal{L}(\xi_0) \in \{\mathcal{H}(\vartheta), \Lambda_\vartheta\}$. As $N \to \infty$, the processes $X^N$ defined in (1.4) converge in finite dimensional distributions to a Gaussian process $X$ (whose covariance structure depends on $d$, $\vartheta$ and the choice of the initial condition, as specified in Theorem 1.1 and in Theorem 1.2).*

The rest of this section is devoted to the proof of Proposition 3.1 in the various cases.

3.1. *The case $d \geq 4$.* Our strategy is as follows: similarly to the technique applied in [24], we are looking for a function $G(\xi)$ that satisfies $LG(\xi) = (\xi(0) - \vartheta) +$ "small error" in order to obtain a representation of the form centered occupation time = martingale + "small error term".

We then use a general functional central limit theorem to treat the martingale term, while we use second moment estimates to show that the error term becomes small. Put

$$G_\lambda(\xi) = \sum_{x \in \mathbb{Z}^d} g_\lambda(x,0)(\xi(x) - \vartheta),$$

where $g_\lambda$ is the resolvent of the underlying random walk. By Lemma 2.2, we have

(3.1) $$(\lambda \mathrm{Id} - \mathcal{L})G_\lambda(\xi) = \xi(0) - \vartheta.$$

Again by Lemma 2.2,

(3.2) $$M_t^\lambda := G_\lambda(\xi_t) - G_\lambda(\xi_0) - \int_0^t \mathcal{L}G_\lambda(\xi_s)\,ds$$

(3.3) $$= \sum_{x,y \in \mathbb{Z}^d} (g_\lambda(y,0) - g_\lambda(x,0))\tilde{N}_t^{x,y} + \sum_{x \in \mathbb{Z}^d} g_\lambda(x,0)(\tilde{N}_t^{x,+} - \tilde{N}_t^{x,-})$$

is a martingale. Using (3.1), we obtain a representation

(3.4) $$\int_0^t (\xi_s(0) - \vartheta)\,ds = -G_\lambda(\xi_t) + G_\lambda(\xi_0) + \lambda \int_0^t G_\lambda(\xi_s)\,ds + M_t^\lambda$$
$$=: R_t^\lambda + M_t^\lambda.$$

We choose $\lambda = 1/N$ and we study the terms $h_d(N)^{-1}R_{Nt}^{1/N}$ and $h_d(N)^{-1}M_{Nt}^{1/N}$ separately in two steps.

*Martingale part*: Using Lemma 2.1, we have

(3.5) $$\langle M^{1/N} \rangle_t = \sum_{x,y \in \mathbb{Z}^d} (g_{1/N}(y,0) - g_{1/N}(x,0))^2 \int_0^t a(x,y)\xi_s(x)\,ds$$
$$+ \sum_{x \in \mathbb{Z}^d} g_{1/N}(x,0)^2 \int_0^t \sigma(\xi_s(x))\,ds.$$

*Case* 1: $(d > 4)$ The martingale $N^{-1/2}M_{Nt}^{1/N}$ has globally bounded jumps $(g_\lambda(x,0) \leq g(x,0) \leq \|g\|_\infty < \infty)$, furthermore, the jump size tends to 0 as $N \to \infty$. (3.5) yields, for any fixed $t > 0$,

$$\langle N^{-1/2}M_{N\cdot}^{1/N} \rangle_t = \sum_{x,y \in \mathbb{Z}^d} a(x,y)(g_{1/N}(y,0) - g_{1/N}(x,0))^2 \frac{1}{N}\int_0^{Nt} \xi_s(x)\,ds$$

$$+ \sum_{x \in \mathbb{Z}^d} g_{1/N}(x,0)^2 \frac{1}{N} \int_0^{Nt} \sigma(\xi_s(x))\, ds$$

$$\xrightarrow[N \to \infty]{P} \text{const} \cdot t.$$

This can be seen as follows. We decompose in both terms the sum over $x \in \mathbb{Z}^d$ into the sum over a ball with a large but fixed radius and the sum over the complement of this ball. For each point $x$ inside the ball, we use that $(1/N) \int_0^{Nt} \xi_s(x)\, ds$ and $(1/N) \int_0^{Nt} \sigma(\xi_s(x))\, ds$ converge to $\vartheta t$ respectively $\sigma_\vartheta^{\text{eq}} t$ in $L_1$ by Proposition 2.3. For $x$ outside the ball, we estimate $\mathbb{E}[\sigma(\xi_s(x))] \leq c_2 \vartheta$ and we use that $g$ is in $\ell_2(\mathbb{Z}^d)$ for $d > 4$. Then the sum over the complement of the ball is small if the radius is large enough. This proves that the r.h.s. converges in $L^1$, so, in particular, it converges in probability.

Using Proposition II.1 in [25], we conclude that $(N^{-1/2} M_{Nt}^{1/N})_{t \geq 0}$ converges in distribution to the law of a Brownian motion.

*Case* 2: ($d = 4$) Here we have to slightly modify our approach because the Green's function is no longer in $\ell_2(\mathbb{Z}^4)$. Instead, we note that

$$\frac{1}{\log N} \sum_{x \in \mathbb{Z}^d} g_{1/N}(x,0)^2 \xrightarrow[N \to \infty]{} \text{const.} > 0$$

and that

$$\frac{1}{\log N} \sum_{x,y \in \mathbb{Z}^d} a(x,y)(g_{1/N}(y,0) - g_{1/N}(x,0))^2$$

$$= \frac{2}{\log N} \langle g_{1/N}(\cdot,0), (-L) g_{1/N}(\cdot,0) \rangle$$

$$= \frac{2}{\log N} \left\langle g_{1/N}(\cdot,0), \delta_0 - \frac{1}{N} g_{1/N}(\cdot,0) \right\rangle$$

$$\leq \frac{2}{\log N} g_{1/N}(0,0) \leq \frac{2}{\log N} g(0,0) \xrightarrow[N \to \infty]{} 0.$$

We then argue analogously to the case above that $(\frac{1}{\sqrt{N \log N}} M_{Nt}^{1/N})_{t \geq 0}$ converges in distribution to a Brownian motion.

*Error part*: Let us first consider $\Lambda_\vartheta$ as the initial condition. We estimate $\mathbb{E}^{\Lambda_\vartheta}[(G_\lambda(\xi_0))^2]$ in order to treat the remainder term. By Lemma 3.3, we have

$$\mathbb{E}^{\Lambda_\vartheta}[(G_\lambda(\xi_0))^2]$$

$$= \sum_{x,y \in \mathbb{Z}^d} g_\lambda(x,0) g_\lambda(y,0) \operatorname{Cov}^{\Lambda_\vartheta}(\xi_0(x), \xi_0(y))$$

$$= \vartheta \int_0^\infty dt\, e^{-\lambda t} \int_0^\infty ds\, e^{-\lambda s} \sum_{x \in \mathbb{Z}^d} a_t(x,0) a_s(x,0)$$

(3.6)
$$+ \frac{\sigma_\vartheta^{\mathrm{eq}}}{2} \int_0^\infty dt\, e^{-\lambda t} \int_0^\infty ds\, e^{-\lambda s} \int_0^\infty du \sum_{x,y \in \mathbb{Z}^d} a_t(x,0) a_s(y,0) a_u(x,y)$$

$$= \int_0^\infty dt \int_0^\infty ds\, e^{-\lambda(t+s)} \left\{ \vartheta a_{t+s}(0,0) + \frac{\sigma_\vartheta^{\mathrm{eq}}}{2} \int_0^\infty du\, a_{t+s+u}(0,0) \right\}$$

$$= \int_0^\infty dr\, e^{-\lambda r} r \left\{ \vartheta a_r(0,0) + \frac{\sigma_\vartheta^{\mathrm{eq}}}{2} \int_r^\infty dv\, a_v(0,0) \right\}.$$

For $d > 4$, we estimate, using (2.2)

$$\mathbb{E}^{\Lambda_\vartheta}[(G_\lambda(\xi_0))^2] \leq C\left(1 + \int_1^\infty dr\, e^{-\lambda r} r(r^{-d/2} + r^{-d/2+1})\right)$$

$$\leq 2C\left(1 + \int_1^\infty dr\, e^{-\lambda r} r^{-d/2+2}\right)$$

to find that

$$\mathbb{E}^{\Lambda_\vartheta}[(N^{-1/2} G_{1/N}(\xi_0))^2] \leq \frac{C}{N} + C \int_1^\infty e^{-r/N} r^{-d/2+2} \frac{dr}{N} = \frac{C}{N} + \frac{C'}{N^{d/2-2}} \xrightarrow[N \to \infty]{} 0.$$

The case $d = 4$ can be treated analogously:

$$\mathbb{E}^{\Lambda_\vartheta}[((N \log N)^{-1/2} G_{1/N}(\xi_0))^2] \leq \frac{C}{N \log N} + \frac{C}{\log N} \int_1^\infty e^{-s/N} \frac{ds}{N} \xrightarrow[N \to \infty]{} 0.$$

Thus, the second term of $R_t$ in (3.4) converges to 0 in $L^2$ after norming with $h_d(N)$, so, in particular, it converges to 0 in probability. By the time-stationarity of $(\xi_t)$ started from $\Lambda_\vartheta$, we see that also the normed first term in (3.4) converges to 0 in probability. Finally, the remaining integral term can be estimated in the following way:

$$\mathbb{E}^{\Lambda_\vartheta} \left| \frac{1}{h_d(N)} \frac{1}{N} \int_0^{Nt} G_{1/N}(\xi_s)\, ds \right| \leq t \mathbb{E}^{\Lambda_\vartheta} |h_d(N)^{-1} G_{1/N}(\xi_0)| \xrightarrow[N \to \infty]{} 0.$$

Putting things together, we conclude that $(h_d(N)^{-1} R_{Nt}^{1/N})_{t \geq 0} \to 0$ as $N \to \infty$ in the sense of finite-dimensional distributions.

Now consider Poisson initial conditions. Note that

$$\mathbb{E}^{\mathcal{H}(\vartheta)}[(G_\lambda(\xi_0))^2] = \vartheta \sum_{x \in \mathbb{Z}^d} g_\lambda(x,0)^2 = \vartheta \int_0^\infty ds \int_0^\infty du\, e^{-\lambda(s+u)} a_{s+u}(0,0),$$

a term which already appeared in (3.6). For the first term in (3.4), note that, by Lemma 3.3,

$$\mathbb{E}^{\mathcal{H}(\vartheta)}[(G_\lambda(\xi_t))^2]$$
$$= \sum_{x,y \in \mathbb{Z}^d} g_\lambda(x,0) g_\lambda(y,0) \operatorname{Cov}^{\mathcal{H}(\vartheta)}(\xi_t(x) \xi_t(y))$$

$$= \vartheta \int_0^\infty ds \int_0^\infty du \, e^{-\lambda(s+u)} a_{s+u}(0,0)$$
$$+ \int_0^\infty ds \int_0^\infty du \int_0^t dr \, e^{-\lambda(s+u)} a_{s+u+2r}(0,0) \mathbb{E}^{\mathcal{H}(\vartheta)}[\sigma(\xi_{t-r}(0))].$$

Estimating $\mathbb{E}^{\mathcal{H}(\vartheta)}[\sigma(\xi_{t-r}(0))] \leq c_2 \vartheta$, we obtain again a term which already appeared in (3.6). Since this estimate is uniform in $t$, we get immediately convergence of the last term in the definition of $R_{Nt}^{1/N}$ in (3.4) to zero.

3.2. *The case $d = 3$.* The decomposition (3.4) of the occupation time in a martingale term and a remainder term as for the case $d > 3$ can not be used in the case $d = 3$: First, $N^{-3/4} G_{1/N}(\xi)$ does not become small in $L^2$, second, as the limit process cannot be a Brownian motion, the Rebolledo-type arguments we used above would not help anyway.

Our approach, again inspired by [24], is to instead "distill" a white noise out of the space–time fluctuations of the ergodic branching random walk system, and to express the normalized occupation time process as a linear functional of this approximate white noise. Technically, for a (momentarily fixed) time horizon $T$, we decompose the occupation time in a term $M_T^T$ and a remainder term, where $M_T^T$ is the final value of a martingale $(M_t^T)_{t \leq T}$.

Recall $u_t(x,0) = \int_0^t a_s(x,0) \, ds$ and define

$$U_t^T(\xi) = \sum_{x \in \mathbb{Z}^3} u_{T-t}(x,0)(\xi(x) - \vartheta).$$

Now

(3.7) $$M_t^T := U_t^T(\xi_t) - U_0^T(\xi_0) - \int_0^t (\partial_s + \mathcal{L}) U_s^T(\xi_s) \, ds$$

is a martingale, and as $(\partial_t - L) u_t(\cdot, 0) = \delta(\cdot, 0)$, we obtain, using Lemma 2.2, the following decomposition of the occupation time:

$$\int_0^T (\xi_s(0) - \vartheta) \, ds = M_T^T + U_0^T(\xi_0).$$

Being interested in $N^{-3/4} \int_0^{NT} (\xi_s(0) - \vartheta) \, ds$, we find ourselves obliged to study $N^{-3/4} M_{NT}^{NT}$ and $N^{-3/4} U_0^{NT}(\xi_0)$.

LEMMA 3.2. *Let $\mathcal{L}(\xi_0) \in \{\Lambda_\vartheta, \mathcal{H}(\vartheta)\}$ for the case of independent branching or $\mathcal{L}(\xi_0) = \mathcal{H}(\vartheta)$ for the case of state dependent branching. The processes*

$$(N^{-3/4} M_{NT}^{NT})_{T \geq 0} \quad and \quad (N^{-3/4} U_0^{NT}(\xi_0))_{T \geq 0}$$

*converge jointly in the sense of finite-dimensional distributions to independent Gaussian limits.*

PROOF. We first consider $(N^{-3/4}U_0^{NT}(\xi_0))_{T\geq 0}$. If we start from a Poisson field, that is, $\mathcal{L}(\xi_0) = \mathcal{H}(\vartheta)$, $N^{-3/4}U^{NT}$ will converge in finite-dimensional distributions to the zero process: The norming with $N^{-3/4}$ is too strong in this case, as can be seen, for example, from

$$N^{-3/2}\mathbb{E}^{\mathcal{H}(\vartheta)}[(U_0^{NT}(\xi_0))^2] = N^{-3/2}\vartheta \sum_{x\in\mathbb{Z}^3} u_{NT}(x,0)^2 = O(N^{-1}).$$

On the other hand, if $\mathcal{L}(\xi_0) = \Lambda_\vartheta$, the norming will be adequate. In the independent branching case the processes $N^{-3/4}U^{NT}$ will have a nontrivial Gaussian limit. Heuristically, if we could simply replace $a_t(0,x)$ by its local CLT analogue, we would find

$$N^{-3/4}U_0^{NT}(\xi_0)$$
$$= N^{-3/4} \sum_{x\in\mathbb{Z}^3} [\xi_0(x) - \vartheta] \int_0^{NT} a_s(0,x)\,ds$$

(3.8)
$$\approx N^{-3/4} \sum_{x\in\mathbb{Z}^3} [\xi_0(x) - \vartheta]$$
$$\times \int_0^{NT} (2\pi s)^{-3/2} (\det Q)^{-1/2} \exp\left(-\frac{x^T Q^{-1} x}{2s}\right) ds$$
$$= N^{-5/4} \sum_{x\in\mathbb{Z}^3} [\xi_0(x) - \vartheta]\varphi_T(x/\sqrt{N}),$$

where $\varphi_T(x) = \int_0^T (2\pi r)^{-3/2}(\det Q)^{-1/2} \exp(-\frac{x^T Q^{-1}x}{2r})\,dr$. If, furthermore, $\varphi_T$ were a Schwartz function, we could conclude using Theorem 1 in [28]. The method of proof used there can be adapted to our situation; technical details can be found in [11], Lemma B.1.

Now let us consider $M_{NT}^{NT}$. Using Lemma 2.2, we can write [we abbreviate $u_s(x) := u_s(x,0)$]

$$M_t^T = \sum_{x,y\in\mathbb{Z}^3} \int_0^t (u_{T-s}(y) - u_{T-s}(x))\,d\tilde{N}_s^{x,y}$$
$$+ \sum_{x\in\mathbb{Z}^3} \int_0^t u_{T-s}(x)\,d\tilde{N}_s^{x,+} - \sum_{x\in\mathbb{Z}^3} \int_0^t u_{T-s}(x)\,d\tilde{N}_s^{x,-}.$$

Now we replace $t$ and $T$ by $NT$ and multiply by $N^{-3/4}$ which yields

$$N^{-3/4}M_{NT}^{NT} = Z_1(N,T) + Z_2(N,T) - Z_3(N,T),$$

where

$$Z_1(N,T) = N^{-3/4} \sum_{x,y\in\mathbb{Z}^3} \int_0^{NT} (u_{NT-s}(y) - u_{NT-s}(x))\,d\tilde{N}_s^{x,y},$$

$$Z_2(N,T) = N^{-3/4} \sum_{x \in \mathbb{Z}^3} \int_0^{NT} u_{NT-s}(x) \, d\tilde{N}_s^{x,+},$$

$$Z_3(N,T) = N^{-3/4} \sum_{x \in \mathbb{Z}^3} \int_0^{NT} u_{NT-s}(x) \, d\tilde{N}_s^{x,-}.$$

We proceed in two steps. In the first step we investigate $Z_1(N,T)$ and in the second step we consider $Z_2(N,T)$ and $Z_3(N,T)$.

*Step* 1: The term $Z_1(N,T)$ converges to zero in probability, since the second moment converges to zero:

$$\mathbb{E}[(Z_1(N,T))^2] = \vartheta N^{-3/2} \sum_{x,y \in \mathbb{Z}^3} a(x,y) \int_0^{NT} (u_{NT-s}(y) - u_{NT-s}(x))^2 \, ds$$

$$= \vartheta N^{-3/2} \sum_{x,y \in \mathbb{Z}^3} a(x,y) \int_0^{NT} (u_s(y) - u_s(x))^2 \, ds$$

$$= 2\vartheta N^{-3/2} \int_0^{NT} \langle u_s, (-Lu_s) \rangle \, ds$$

$$= 2\vartheta N^{-3/2} \int_0^{NT} \langle u_s, \delta_0 - a_s(\cdot, 0) \rangle \, ds$$

$$\leq 2\vartheta N^{-3/2} \int_0^{NT} u_s(0) \, ds \leq \vartheta N^{-1/2} T g(0,0) \to 0.$$

*Step* 2: Now we consider the remaining terms $Z_2(N,T)$ and $Z_3(N,T)$. We define a random field $Y_{N,T}$ on $L^2([0,T] \times \mathbb{R}^3)$ via

$$\langle Y_{N,T}, \varphi \rangle := N^{1/4} \int_{\mathbb{R}^3} dz \int_0^T d\tilde{N}_{Ns}^{\sqrt{N}\lfloor z \rfloor_N, +} \varphi(s,z)$$

$$= N^{1/4} \sum_{x \in \mathbb{Z}^3/\sqrt{N}} \int_0^T d\tilde{N}_{Ns}^{\sqrt{N}x,+} \int_{x+\Omega_N} dz \, \varphi(s,z),$$

where $\lfloor z \rfloor_N$ is determined by $\lfloor z \rfloor_N \in \mathbb{Z}^3/\sqrt{N}$ and $z \in \lfloor z \rfloor_N + \Omega_N$, with $\Omega_N = (-\frac{1}{2\sqrt{N}}, \frac{1}{2\sqrt{N}}]^3$. Thus, we can write

$$Z_2(N,T) = \langle Y_{N,T}, v_{N,T} \rangle,$$

where

$$v_{N,T}(s,z) = N^{1/2} \sum_{x \in \mathbb{Z}^3/\sqrt{N}} u_{N(T-s)}(\sqrt{N}x) \mathbf{1}_{x+\Omega_N}(z).$$

Next we wish to show that $Y_{N,T}$ converges toward a white noise $Y_T$ on $[0,T] \times \mathbb{R}^3$ (with covariance measure given by $\sigma_\vartheta^{\text{eq}}/2$ times the Lebesgue measure). Furthermore, for large $N$, the CLT suggests that $v_{N,T}$ should be similar to

$$v_T(s,z) := \int_0^{T-s} p_r(z,0)\,dr,$$

where $p_r(x,y) := (2\pi r)^{-3/2}(\det Q)^{-1/2}\exp(-\frac{(y-x)^T Q^{-1}(y-x)}{2r})$. Thus, we expect $Z_2(N,T) \approx \langle Y_T, v_T \rangle$, which shows the Gaussian nature. We proceed in two parts to justify these heuristics:

*Part* 1: Here we show that $\langle Y_{N,T}, \varphi \rangle \to \langle Y_T, \varphi \rangle$ as $N \to \infty$ when $\varphi \in L^2_{(1/2)\sigma^{\text{eq}}_\vartheta}([0,T] \times \mathbb{R}^3)$. The index $\frac{1}{2}\sigma^{\text{eq}}_\vartheta$ indicates that this is the $L^2$-space corresponding to $\frac{1}{2}\sigma^{\text{eq}}_\vartheta$ times the Lebesgue measure on $[0,T] \times \mathbb{R}^3$. We write $\|\varphi\|_2$ for the norm of $\varphi$ in this space. $Y_T$ is a space–time white noise based on $\frac{1}{2}\sigma^{\text{eq}}_\vartheta$ times the Lebesgue measure on $[0,T] \times \mathbb{R}^3$, that is, a $Y_T : \varphi \mapsto \langle \varphi, Y_T \rangle$ is a linear isometry from $L^2_{(1/2)\sigma^{\text{eq}}_\vartheta}([0,T] \times \mathbb{R}^3)$ to the space of Gaussian random variables equipped with the $L^2$-norm. See, for example, Chapter 1 of [27] for background on white noises.

First we consider test functions consisting only of finitely many steps: Let

$$(3.9) \quad \varphi(s,x) = \sum_{k=1}^n \mathbf{1}_{\bigcup_{l=1}^{m(k)}[r_l^k,t_l^k]}(s)\mathbf{1}_{A_k}(x) = \sum_{k=1}^n \mathbf{1}_{A_k}(x)\sum_{l=1}^{m(k)}\mathbf{1}_{[r_l^k,t_l^k]}(s),$$

where $A_1, \ldots, A_n \subset \mathbb{R}^3$ are disjoint (say, bounded closed parallelepipeds) and $r_1^k < t_1^k \leq r_2^k < t_2^k \leq \cdots \leq r_{m(k)}^k < t_{m(k)}^k$. Let

$$Z_t^{N,k} = N^{1/4}\sum_{x \in \mathbb{Z}^3/\sqrt{N}} \lambda(A_k \cap (x + \Omega_N))\tilde{N}_{Nt}^{\sqrt{N}x,+}, \qquad k = 1, \ldots, n.$$

Then $(Z_t^N)_{0 \leq t \leq T} = (Z_t^{N,1}, \ldots, Z_t^{N,n})_{0 \leq t \leq T}$ is an $\mathbb{R}^n$-valued martingale. The assumptions of Proposition II.1 in [25] are fulfilled since:

(i) We observe for $k \neq l$ that

$$\langle Z^{N,k}, Z^{N,l}\rangle_t = N^{1/2}\sum_{x \in \mathbb{Z}^3/\sqrt{N}}\lambda(A_k \cap (x+\Omega_N))\lambda(A_l \cap (x+\Omega_N))$$
$$\times \int_0^{Nt} \tfrac{1}{2}\sigma(\xi_s(\sqrt{N}x))\,ds,$$

which converges in probability to 0 since

$$\mathbb{E}[\langle Z^{N,k}, Z^{N,l}\rangle_t] \leq \tfrac{1}{2}c_2 \vartheta t N^{-3/2}$$
$$\times \#\{x \in \mathbb{Z}^3/\sqrt{N} : A_k \cap (x+\Omega_N) \neq \varnothing, A_l \cap (x+\Omega_N) \neq \varnothing\},$$

which is 0 if $N$ is large enough since $A_k$ and $A_l$ are closed and disjoint.

For $k = l$, we calculate

$$\langle Z^{N,k}, Z^{N,k} \rangle_t = N^{1/2} \sum_{x \in \mathbb{Z}^3/\sqrt{N}} (\lambda(A_k \cap (x + \Omega_N)))^2 \int_0^{Nt} \tfrac{1}{2}\sigma(\xi_s(\sqrt{N}x))\, ds$$

$$\xrightarrow[N \to \infty]{P} \tfrac{1}{2} \sigma_\vartheta^{\mathrm{eq}} \lambda(A_k) t,$$

which can be seen by the following argument:

$$\mathbb{E}\left[\left(N^{1/2} \sum_{x \in \mathbb{Z}^3/\sqrt{N}} (\lambda(A_k \cap (x + \Omega_N)))^2 \int_0^{Nt} \sigma(\xi_s(\sqrt{N}x))\, ds - \sigma_\vartheta^{\mathrm{eq}} \lambda(A_k) t \right)^2\right]$$

$$(3.10) \quad = \sum_{x,y \in (\mathbb{Z}^3/\sqrt{N}) \cap A_k} N^{-5} \int_0^{Nt} ds \int_0^{Nt} du$$

$$\times \mathbb{E}[(\sigma(\xi_s(\sqrt{N}x)) - \sigma_\vartheta^{\mathrm{eq}})(\sigma(\xi_u(\sqrt{N}y)) - \sigma_\vartheta^{\mathrm{eq}})] + r_N,$$

where $r_N \to 0$ as $N \to \infty$. Proposition 1 of [9] yields

$$\sup_{z \in \mathbb{Z}^3} \left| \mathbb{E}^{\mathcal{H}(\vartheta)}[f(\xi_s)g(\xi_{s+u}(z + \cdot))] - \int f\, d\Lambda_\vartheta \int g\, d\Lambda_\vartheta \right| \longrightarrow 0 \quad \text{as } s, u \to \infty.$$

From this and shift-invariance of $\xi_t$ we can conclude that the expression on the r.h.s. of (3.10) converges to zero.

(ii) We observe that $Z^{N,k}$ has jumps of size $N^{-5/4}$, such that condition (ii) of Proposition II.1 in [25] is obviously fulfilled.

By Proposition II.1 in [25], we can conclude

$$(3.11) \qquad (Z^{N,1}, \ldots, Z^{N,n}) \xrightarrow[N \to \infty]{} (Z^1, \ldots, Z^n),$$

where $\{Z^k\}$ are independent Brownian motions with variance parameter $\tfrac{1}{2}\sigma_\vartheta^{\mathrm{eq}}\lambda(A_k)$.

For $\varphi$ defined in (3.9), we obtain

$$\langle Y_{N,T}, \varphi \rangle = \sum_{k=1}^n \sum_{l=1}^{m(k)} (Z^{N,k}_{t_l^k} - Z^{N,k}_{r_l^k}) \xrightarrow[N \to \infty]{} \sum_{k=1}^n \sum_{l=1}^{m(k)} (Z^k_{t_l^k} - Z^k_{r_l^k}).$$

The limit is a sum of independent normal random variables by (3.11). Therefore, the limit is normal with variance $\sum_{k=1}^n \sum_{l=1}^{m(k)} \tfrac{1}{2}\sigma_\vartheta^{\mathrm{eq}}\lambda(A_k)(t_l^k - r_l^k)$ and, hence,

$$\langle Y_{N,T}, \varphi \rangle \xrightarrow[N \to \infty]{} \langle Y_T, \varphi \rangle.$$

Then we can extend the convergence statement to all $\varphi \in L^2_{(1/2)\sigma_\vartheta^{\mathrm{eq}}}([0,T] \times \mathbb{R}^3)$, since the functions of the form (3.9) are dense in $L^2_{(1/2)\sigma_\vartheta^{\mathrm{eq}}}([0,T] \times \mathbb{R}^3)$ and since (for $\varphi\colon\mathbb{R}^3 \to \mathbb{R}$ such that $\varphi^2$ is Riemann integrable)

$$\lim_{N\to\infty} \mathbb{E}[\langle Y_{N,T}, \varphi\rangle^2] = \tfrac{1}{2}\sigma^{\mathrm{eq}}_\vartheta \int_{\mathbb{R}^3} \int_0^T \varphi(s,z)^2\, ds\, dz = \|\varphi\|_2^2.$$

The last assertion can be seen by the following calculation (note that $\mathbb{E}[\sigma(\xi_{Ns}(\sqrt{N}x))]$ does not depend on $x$):

$\mathbb{E}[\langle Y_{N,T}, \varphi\rangle^2]$

$$= N^{1/2} \sum_{x\in\mathbb{Z}^3/\sqrt{N}} \int_0^T ds \left[\int_{x+\Omega_N} dz\, \varphi(s,z)\right]^2 N\tfrac{1}{2}\mathbb{E}[\sigma(\xi_{Ns}(\sqrt{N}x))]$$

$$= N^{1/2} \sum_{x\in\mathbb{Z}^3/\sqrt{N}} \int_0^T ds \left[\int_{x+\Omega_N} dz\, \varphi(s,z)\right]^2 N\tfrac{1}{2}\sigma^{\mathrm{eq}}_\vartheta + r_N \qquad \text{(with } r_N \to 0\text{)}$$

$$= \tfrac{1}{2}\sigma^{\mathrm{eq}}_\vartheta \sum_{x\in\mathbb{Z}^3/\sqrt{N}} \int_0^T ds \left[\int_{x+\Omega_N} dz\, \varphi(s,z)\right] \left[N^{3/2}\int_{x+\Omega_N} dw\, \varphi(s,w)\right] + r_N,$$

which is a Riemannian sum for $\tfrac{1}{2}\sigma^{\mathrm{eq}}_\vartheta \int_0^T \int_{\mathbb{R}^3} \varphi^2(s,z)\, dz\, ds$. This completes the proof of the assertion.

*Part 2*: Now we show that $\langle Y_{N,T}, v_{N,T} - v_T\rangle \xrightarrow[N\to\infty]{P} 0$, in fact we will show that

(3.12) $$\mathbb{E}[\langle Y_{N,T}, v_{N,T} - v_T\rangle^2] \xrightarrow[N\to\infty]{} 0.$$

We have

$\mathbb{E}[\langle Y_{N,T}, v_{N,T} - v_T\rangle^2]$

$$= N^{3/2} \sum_{x\in\mathbb{Z}^3/\sqrt{N}} \int_0^T ds \left(\int_{x+\Omega_N} dz\, (v_{N,T}(s,z) - v_T(s,z))\right)^2$$

$$\times \tfrac{1}{2}\mathbb{E}[\sigma(\xi_{Ns}(\sqrt{N}x))]$$

$$\leq \tfrac{1}{2} c_2 \vartheta N^{3/2}$$

$$\times \sum_{x\in\mathbb{Z}^3/\sqrt{N}} \int_0^T ds \left(\int_0^{T-s} dr \int_{x+\Omega_N} dz\, (N^{3/2} a_{Nr}(\sqrt{N}x, 0) - p_r(z,0))\right)^2$$

$$= \tfrac{1}{2} c_2 \vartheta N^{3/2} \sum_{x\in\mathbb{Z}^3/\sqrt{N}} \int_0^T ds \left(\int_0^\varepsilon dr \ldots + \int_\varepsilon^{T-s} dr \ldots\right)^2$$

$$\leq c_2 \vartheta N^{3/2} \sum_{x\in\mathbb{Z}^3/\sqrt{N}} \int_0^T ds \left\{\left[\int_0^\varepsilon dr \ldots\right]^2 + \left[\int_\varepsilon^{T-s} dr \ldots\right]^2\right\}.$$

Now note that

$$N^{3/2} \sum_{x\in\mathbb{Z}^3/\sqrt{N}} \int_0^T ds \left[\int_0^\varepsilon dr \int_{x+\Omega_N} dz\, N^{3/2} a_{Nr}(\sqrt{N}x,0)\right]^2$$

$$= N^{3/2} T \sum_{x\in\mathbb{Z}^3/\sqrt{N}} \left[\frac{1}{N}\int_0^{N\varepsilon} dr\, a_r(\sqrt{N}x,0)\right]^2$$

$$= N^{-1/2} T \|u_{N\varepsilon}(\cdot,0)\|_2^2 \leq C\sqrt{\varepsilon},$$

where we use for the last estimate that

$$\|u_t(\cdot,0)\|_2^2 = \sum_{x\in\mathbb{Z}^d} \int_0^t dr \int_0^t ds\, a_r(x,0) a_s(x,0)$$

$$= 2\int_0^t dr \int_r^t ds\, a_{r+s}(0,0) \sim \text{const.} \times \sqrt{t}$$

by (2.2). Similarly,

$$N^{3/2} \sum_{x\in\mathbb{Z}^3/\sqrt{N}} \int_0^T ds \left[\int_0^\varepsilon dr \int_{x+\Omega_N} dz\, p_r(z,0)\right]^2$$

$$= N^{3/2} T \sum_{x\in\mathbb{Z}^3/\sqrt{N}} \left[\int_0^\varepsilon dr \int_{x+\Omega_N} dz\, p_r(z,0)\right]^2$$

$$\leq N^{3/2} T \sum_{x\in\mathbb{Z}^3/\sqrt{N}} |\Omega_N| \int_{x+\Omega_N} dz \left(\int_0^\varepsilon dr\, p_r(z,0)\right)^2$$

$$= T \int_{\mathbb{R}^3} dz \int_0^\varepsilon dr \int_0^\varepsilon ds\, p_r(z,0) p_s(z,0)$$

$$= 2T \int_0^\varepsilon dr \int_r^\varepsilon ds\, p_{r+s}(0,0) \leq C\sqrt{\varepsilon},$$

where we used the Cauchy–Schwarz inequality.

In order to treat the remaining term we use that [see, e.g. (2.1)]

(3.13) $\quad |N^{3/2} a_{Nr}(\sqrt{N}x,0) - p_r(z,0)| \leq C_\varepsilon \frac{1}{1+|z|^2/r}\psi(N)$

uniformly in $N$, $r \in [\varepsilon, T]$, $x \in \mathbb{Z}^3/\sqrt{N}$, $z \in x + \Omega_N$, where $\psi(N) \to 0$ as $N \to \infty$. (Note that this requires only a second moment assumption on $a$.) This yields

$$N^{3/2} \sum_{x\in\mathbb{Z}^3/\sqrt{N}} \int_0^T ds \left[\int_\varepsilon^{T-s} dr \int_{x+\Omega_N} dz (N^{3/2} a_{Nr}(\sqrt{N}x,0) - p_r(z,0))\right]^2$$

$$\leq N^{3/2}\psi(N)^2 T \sum_{x\in\mathbb{Z}^3/\sqrt{N}} \left(\int_\varepsilon^T dr \int_{x+\Omega_N} dz \frac{C_\varepsilon}{1+|z|^2/r}\right)^2$$

$$\leq N^{3/2}\psi(N)^2 T \sum_{x\in\mathbb{Z}^3/\sqrt{N}} \left(\int_{x+\Omega_N} dz \frac{C_\varepsilon T}{1+|z|^2/T}\right)^2$$

$$\leq N^{3/2}\psi(N)^2 T^3 C_\varepsilon^2 \sum_{x\in\mathbb{Z}^3/\sqrt{N}} |\Omega_N| \int_{x+\Omega_N} dz\, (1+|z|^2/T)^{-2}$$

$$= C_\varepsilon^2 T^3 \psi(N)^2 \int_{\mathbb{R}^3} dz\, (1+|z|^2/T)^{-2} \xrightarrow[N\to\infty]{} 0.$$

Combining, we see that $\limsup_{N\to\infty} \mathbb{E}[\langle Y_{N,T}, v_{N,T} - v_T\rangle^2] \leq C\sqrt{\varepsilon}$, now let $\varepsilon \to 0$ to obtain (3.12).

Thus, we have shown that $Z_2(N,T)$ converges to a Gaussian limit. $Z_3(N,T)$ can be treated completely analogously, and as it involves only integrals with respect to $(\tilde{N}_t^{x,-})$, $x \in \mathbb{Z}^3$, and the martingales $\tilde{N}^{x,-}$ and $\tilde{N}^{x,+}$ are all pairwise orthogonal, we see that $Z_2(N,T)$ and $Z_3(N,T)$ converge jointly to (independent) Gaussian processes. Thus, $(N^{-3/4}M_{NT}^{NT})$ converges as $N \to \infty$ to a Gaussian process.

Finally, a remark on the *joint* convergence of $U^{NT}$ and $M^{NT}$ when starting from the invariant distribution $\Lambda_\vartheta$ is in order: Note that $U_{NT}^{NT}(\xi_0)$ depends only on the initial condition, whereas $M^{NT}$ is a function of the driving martingales $\tilde{N}^{x,\pm}$, $x \in \mathbb{Z}^3$. Scrutinizing the proof the reader will find that even conditional on $\xi_0 = \eta$, $M^{NT}$ will converge to the same Gaussian process, as long as $\eta$ is such that $\mathcal{L}(\xi_t|\xi_0 = \eta) \Rightarrow \Lambda_\vartheta$ as $t \to \infty$. Note that $\Lambda_\vartheta$-a.a. initial conditions $\xi_0$ have this property because $\Lambda_\vartheta$ is an extremal equilibrium in the classical branching case. The argument can be made precise by considering a joint characteristic functional of $M^{NT}$ and $U^{NT}$ and then conditioning on $\xi_0$. □

3.3. *Covariance computation.* First we need a second moment formula for branching random walks.

LEMMA 3.3. *For $u \leq v$, $x, y \in \mathbb{Z}^d$, we have the following moment formulas:*

$$\mathbb{E}^{\mathcal{H}(\vartheta)}[\xi_u(x)\xi_v(y)] = \vartheta^2 + \vartheta a_{v-u}(x,y)$$
$$+ \int_0^u a_{v-u+2r}(x,y) \mathbb{E}^{\mathcal{H}(\vartheta)}[\sigma(\xi_{u-r}(0))]\, dr,$$

$$\mathbb{E}^{\Lambda_\vartheta}[\xi_u(x)\xi_v(y)] = \vartheta^2 + \vartheta a_{v-u}(x,y) + \tfrac{1}{2}\sigma_\vartheta^{\mathrm{eq}} \int_{v-u}^\infty a_r(x,y)\, dr.$$

PROOF. In [8], Lemma 4 one can find the moment formula, but only for the process in different space points at the same time. To obtain the formula for different times $u < v$, simply condition on the configuration $\xi_u$, then use the Markov property and the first moment formula from Lemma 4 in [8], finally use the second moment formula from the same lemma. Moreover, we use that $\mathbb{E}^{\mathcal{H}(\vartheta)}[\sigma(\xi_{u-r}(x))]$ does not depend on $x$. $\square$

Now we compute the covariance of the limit of the renormalized occupation time.

PROPOSITION 3.4. *The variance of the limit of the renormalized occupation time is*

$$\mathbb{E}^\mu[X_s^N X_t^N]$$
$$\xrightarrow[N\to\infty]{} \begin{cases} \dfrac{\sqrt{2}}{3\pi^{3/2}}(\det Q)^{-1/2}\sigma_\vartheta^{\mathrm{eq}}[t^{3/2} + s^{3/2} - |t-s|^{3/2}], \\ \qquad\qquad\qquad\qquad\qquad d=3, \mu = \Lambda_\vartheta, \\ \dfrac{2\sqrt{2}}{3\pi^{3/2}}(\det Q)^{-1/2}\sigma_\vartheta^{\mathrm{eq}}\left[t^{3/2} + s^{3/2} - \dfrac{1}{2}|t-s|^{3/2} - \dfrac{1}{2}(t+s)^{3/2}\right], \\ \qquad\qquad\qquad\qquad\qquad d=3, \mu = \mathcal{H}(\vartheta), \\ (2\pi)^{-2}(\det Q)^{-1/2}\sigma_\vartheta^{\mathrm{eq}} \times (s \wedge t), \qquad d=4, \mu \in \{\Lambda_\vartheta, \mathcal{H}(\vartheta)\}, \\ \left[2\vartheta \int_0^\infty du\, a_u(0,0) + \sigma_\vartheta^{\mathrm{eq}} \int_0^\infty du\, u a_u(0,0)\right] \times (s \wedge t), \\ \qquad\qquad\qquad\qquad\qquad d\geq 5, \mu \in \{\Lambda_\vartheta, \mathcal{H}(\vartheta)\}. \end{cases}$$

PROOF. The proof is split up into different cases. We assume $s \leq t$ throughout. Let us first consider the situation $\mathcal{L}(\xi_0) = \mathcal{H}(\vartheta)$. By Lemma 3.3, we have

$$\mathbb{E}^{\mathcal{H}(\vartheta)}[X_s^N X_t^N] = \frac{1}{h_d(N)^2} \int_0^{Ns} du \int_0^{Nt} dv\, \mathrm{Cov}^{\mathcal{H}(\vartheta)}(\xi_u(0), \xi_v(0))$$

(3.14)
$$= \frac{\vartheta}{h_d(N)^2} \int_0^{Ns} du \int_0^{Nt} dv\, a_{|v-u|}(0,0)$$
$$+ \frac{1}{h_d(N)^2} \int_0^{Ns} du \int_0^{Nt} dv \int_0^{u \wedge v} dr$$
$$\times a_{|v-u|+2r}(0,0) \mathbb{E}^{\mathcal{H}(\vartheta)}[\sigma(\xi_{u \wedge v - r}(0))]$$
$$=: I_1 + I_2.$$

*Case* 1: Let $d = 3$. We have $0 \leq I_1 \leq \vartheta N^{-3/2}(Ns) \int_0^\infty a_r(0,0)\, dr = O(N^{-1/2})$, so that this term is asymptotically negligible. Fix $\varepsilon > 0$ for the moment. By (2.2), we can find $K > 0$ such that $a_r(0,0) \leq (1+\varepsilon)c_3 r^{-3/2}$ and $\mathbb{E}^{\mathcal{H}(\vartheta)}[\sigma(\xi_r(0))] \leq$ $(1+\varepsilon)\sigma_\vartheta^{\text{eq}}$ for $r \geq K$, where $c_3 = (2\pi)^{-3/2}(\det Q)^{-1/2}$. Thus, we can bound $I_2$ by

$$
\begin{aligned}
&\frac{2(1+\varepsilon)^2\sigma_\vartheta^{\text{eq}}c_3}{N^{3/2}} \int_0^{Ns-K} du \int_{u+K}^{Ns} dv \int_0^u \frac{dr}{(v-u+2r)^{3/2}} \\
&\quad + \frac{(1+\varepsilon)^2\sigma_\vartheta^{\text{eq}}c_3}{N^{3/2}} \int_0^{Ns-K} du \int_{Ns}^{Nt} dv \int_0^u \frac{dr}{(v-u+2r)^{3/2}} \\
&\quad + O(N^{-1/2}).
\end{aligned}
$$ (3.15)

The first term in (3.15) is equal to

$$
\begin{aligned}
&\frac{2(1+\varepsilon)^2\sigma_\vartheta^{\text{eq}}c_3}{N^{3/2}} \int_0^{Ns-K} du \int_{u+K}^{Ns} dv \left\{\frac{1}{(v-u)^{1/2}} - \frac{1}{(v+u)^{1/2}}\right\} \\
&= \frac{4(1+\varepsilon)^2\sigma_\vartheta^{\text{eq}}c_3}{N^{3/2}} \Bigg[-\frac{2}{3}(Ns-u)^{3/2} - K^{1/2}u \\
&\qquad\qquad - \frac{2}{3}(Ns+u)^{3/2} + \frac{1}{3}(2u+K)^{3/2}\Bigg]_{u=0}^{u=Ns-K} \\
&\xrightarrow[N\to\infty]{} \frac{4}{3}(1+\varepsilon)^2\sigma_\vartheta^{\text{eq}}c_3(4 - 2^{3/2})s^{3/2}.
\end{aligned}
$$

Analogously, the second term in (3.15) is equal to

$$
\begin{aligned}
&\frac{2(1+\varepsilon)^2\sigma_\vartheta^{\text{eq}}c_3}{N^{3/2}} \int_0^{Ns-K} du \left\{[(v-u)^{1/2} - (v+u)^{1/2}]_{v=Ns}^{v=Nt}\right\} \\
&= \frac{2(1+\varepsilon)^2\sigma_\vartheta^{\text{eq}}c_3}{N^{3/2}} \Bigg[-\frac{2}{3}(Nt-u)^{3/2} + \frac{2}{3}(Ns-u)^{3/2} \\
&\qquad\qquad - \frac{2}{3}(Nt+u)^{3/2} + \frac{2}{3}(Ns+u)^{3/2}\Bigg]_{u=0}^{u=Ns-K} \\
&\xrightarrow[N\to\infty]{} \frac{4}{3}(1+\varepsilon)^2\sigma_\vartheta^{\text{eq}}c_3(2t^{3/2} - (t-s)^{3/2} - (t+s)^{3/2} - (2 - 2^{3/2})s^{3/2}).
\end{aligned}
$$

Combining these terms and letting $\varepsilon \to 0$, we see that

$$
\limsup_{N\to\infty} I_2 \leq \tfrac{8}{3}c_3\sigma_\vartheta^{\text{eq}}(t^{3/2} + s^{3/2} - \tfrac{1}{2}(t-s)^{3/2} - \tfrac{1}{2}(t+s)^{3/2}).
$$

$\liminf_{N\to\infty} I_2$ can be analogously bounded from below, concluding the proof in this case.

*Case* 2: Let $d = 4$. We have $0 \leq I_1 \leq \vartheta(N\log N)^{-1}(Ns)\int_0^\infty a_r(0,0)\,dr = O(1/\log N)$, so that this term is again asymptotically negligible. Arguing as in case 1, we can now bound $I_2$ from above by

$$
\frac{2(1+\varepsilon)^2\sigma_\vartheta^{\text{eq}}c_4}{N\log N} \int_0^{Ns-K} du \int_{u+K}^{Ns} dv \int_0^u \frac{dr}{(v-u+2r)^2}
$$

$$\text{(3.16)} \qquad + \frac{(1+\varepsilon)^2 \sigma_\vartheta^{\text{eq}} c_4}{N \log N} \int_0^{Ns-K} du \int_{Ns}^{Nt} dv \int_0^u \frac{dr}{(v-u+2r)^2}$$

$$+ O((\log N)^{-1}),$$

where $c_4 = (2\pi)^{-2} (\det Q)^{-1/2}$. The first term in (3.16) is equal to

$$\frac{(1+\varepsilon)^2 \sigma_\vartheta^{\text{eq}} c_4}{N \log N} \int_0^{Ns-K} du \int_{u+K}^{Ns} dv \left( \frac{1}{v-u} - \frac{1}{v+u} \right)$$

$$= \frac{(1+\varepsilon)^2 \sigma_\vartheta^{\text{eq}} c_4}{N \log N} \int_0^{Ns-K} du (\log(Ns-u)$$

$$- \log K - \log(Ns+u) + \log(2u+K))$$

$$= \frac{(1+\varepsilon)^2 \sigma_\vartheta^{\text{eq}} c_4}{N \log N} [-(Ns-u) \log(Ns-u) - u - u \log K$$

$$- (Ns+u) \log(Ns+u)$$

$$+ u + (u+K/2) \log(2u+K) - u]_{u=0}^{u=Ns-K}$$

$$\longrightarrow (1+\varepsilon)^2 \sigma_\vartheta^{\text{eq}} c_4 s \qquad \text{as } N \to \infty.$$

Now the second term in (3.16) is bounded above by

$$\frac{(1+\varepsilon)^2 \sigma_\vartheta^{\text{eq}} c_4}{2N \log N} \int_0^{Ns-K} du \int_{Ns}^{Nt} \frac{dv}{v-u}$$

$$= \frac{(1+\varepsilon)^2 \sigma_\vartheta^{\text{eq}} c_4}{2N \log N} [-(Nt-u) \log(Nt-u) - u$$

$$+ (Ns-u) \log(Ns-u) + u]_{u=0}^{u=Ns-K}$$

$$= \frac{(1+\varepsilon)^2 \sigma_\vartheta^{\text{eq}} c_4}{2N \log N} \{ Nt \log(Nt) - (N(t-s)+K) \log(N(t-s)+K)$$

$$+ K \log K - Ns \log(Ns) \}$$

$$\xrightarrow[N \to \infty]{} 0.$$

Thus, letting $\varepsilon \to 0$, we see that $\limsup_{N \to \infty} I_2 \le \sigma_\vartheta^{\text{eq}} c_4 s$ in this case. Again, $\liminf_{N \to \infty} I_2$ can be bounded analogously, completing this part of the proof.

*Case* 3: Let $d \ge 5$. We have

$$I_1 = \frac{2\vartheta}{N} \int_0^{Ns} du \int_u^{Ns} dv\, a_{v-u}(0,0) + \frac{\vartheta}{N} \int_0^{Ns} du \int_{Ns}^{Nt} dv\, a_{v-u}(0,0)$$

$$= \frac{2\vartheta}{N} \int_0^{Ns} du \int_0^\infty dr\, a_r(0,0) + O\left( N^{-1} \int_0^{Ns} du \int_{Ns}^\infty dv\, a_{v-u}(0,0) \right)$$

$$= 2\vartheta s \int_0^\infty dr\, a_r(0,0) + O(1/N).$$

We decompose $I_2$ as

$$I_2 = \frac{2}{N} \int_0^{Ns} du \int_u^{Ns} dv \int_0^u dr \, a_{v-u+2r}(0,0) \mathbb{E}^{\mathcal{H}(\vartheta)}[\sigma(\xi_{u-r}(0))]$$
$$+ \frac{1}{N} \int_0^{Ns} du \int_{Ns}^{Nt} dv \int_0^u dr \, a_{v-u+2r}(0,0) \mathbb{E}^{\mathcal{H}(\vartheta)}[\sigma(\xi_{u-r}(0))]. \tag{3.17}$$

The second term in (3.17) is bounded by

$$\frac{C}{N} \int_{Ns}^{Nt} dv \int_{v-Ns}^{\infty} dr \, a_r(0,0) + \frac{C}{N} \int_0^{Ns-1} du \int_{Ns}^{Nt} dv \int_{v-u}^{\infty} dr \, a_r(0,0)$$
$$= O(N^{-1}) + O(N^{2-d/2})$$

by (2.2).

Choose $K$ large enough such that $\mathbb{E}^{\mathcal{H}(\vartheta)}[\sigma(\xi_r(0))] \leq (1+\varepsilon)\sigma_\vartheta^{\text{eq}}$ for all $r \geq K$. The first term in (3.17) can be bounded by

$$\frac{2(1+\varepsilon)\sigma_\vartheta^{\text{eq}}}{N} \int_K^{Ns} du \int_u^{Ns} dv \int_0^{u-K} dr \, a_{v-u+2r}(0,0) + O(N^{-1})$$
$$= (1+\varepsilon)\sigma_\vartheta^{\text{eq}} \int_{K/N}^{s} du \int_0^{N(s-u)} dv \int_v^{v+2Nu-2K} dr \, a_r(0,0) + O(N^{-1})$$
$$\longrightarrow (1+\varepsilon)\sigma_\vartheta^{\text{eq}} s \int_0^\infty dv \int_v^\infty dr \, a_r(0,0) = (1+\varepsilon)\sigma_\vartheta^{\text{eq}} s \int_0^\infty du \, u a_u(0,0).$$

Letting $\varepsilon \to 0$, we see that $\limsup_{N\to\infty}$ of the first term in (3.17) is bounded by $\sigma_\vartheta^{\text{eq}} s \int_0^\infty du \, u a_u(0,0)$. Since the liminf can be bounded analogously, we are done with the case $d \geq 5$.

Now let us consider the situation $\mathcal{L}(\xi_0) = \Lambda_\vartheta$. Then we have, by Lemma 3.3,

$$\mathbb{E}^{\Lambda_\vartheta}[X_s^N X_t^N] = \frac{\vartheta}{h_d(N)^2} \int_0^{Ns} du \int_0^{Nt} dv \, a_{|v-u|}(0,0)$$
$$+ \frac{\sigma_\vartheta^{\text{eq}}}{2h_d(N)^2} \int_0^{Ns} du \int_0^{Nt} dv \int_{|v-u|}^\infty dr \, a_r(0,0) \tag{3.18}$$
$$=: I_1 + I_2'.$$

The computations for $d = 4$ and $d \geq 5$ are entirely analogous to those above, and will be omitted. Let us briefly comment on the case $d = 3$ in this situation. $I_1$ is again negligible, and choosing $K$ large enough, we can now bound $I_2'$ from above by

$$\frac{2(1+\varepsilon)\sigma_\vartheta^{\text{eq}} c_3}{N^{3/2}} \int_0^{Ns-K} du \int_{u+K}^{Ns} \frac{dv}{(v-u)^{1/2}}$$

$$+ \frac{(1+\varepsilon)\sigma_\vartheta^{\mathrm{eq}} c_3}{N^{3/2}} \int_0^{Ns-K} du \int_{Ns}^{Nt} \frac{dv}{(v-u)^{1/2}} + O(N^{-1/2})$$

$$= \frac{2(1+\varepsilon)\sigma_\vartheta^{\mathrm{eq}} c_3}{N^{3/2}} \int_0^{Ns-K} du\, [(Ns-u)^{1/2} + (Nt-u)^{1/2}] + O(N^{-1/2})$$

$$\xrightarrow[N\to\infty]{} \frac{4}{3}(1+\varepsilon)\sigma_\vartheta^{\mathrm{eq}} c_3(t^{3/2} + s^{3/2} - (t-s)^{3/2}).$$

We conclude the proof as above. □

**4. Proof of Theorem 1.1.** Here we complete the proof of Theorem 1.1. In view of Proposition 3.1, it suffices to check that the sequence $X^N$, $N \in \mathbb{N}$, is tight (e.g., in the space of all continuous processes, equipped with the norm of locally uniform convergence). In order to do so, we use the well-known criterion on moments of increments, stating that a sequence of processes $X^N$ is tight (and, furthermore, any limit point has continuous paths) if there exist $\alpha, \beta > 0$ such that, for each $t_0 > 0$,

(4.1) $$\mathbb{E}[|X_t^N - X_s^N|^\alpha] \leq C(t-s)^{1+\beta}$$

holds uniformly in $N$ and $0 \leq s < t \leq t_0$ (see, e.g., [18], Corollary 14.9).

For $d \geq 4$ and independent branching, with $\alpha = 4$ and $\beta = 1$, this is the content of the following lemma.

LEMMA 4.1. *Let $d \geq 4$ and $\mu \in \{\mathcal{H}(\vartheta), \Lambda_\vartheta\}$. Then in the independent branching case for each $t_0 > 0$ there is a $C = C(t_0, \vartheta, d)$ such that, for all $0 \leq s \leq t \leq t_0$, we have*

$$\mathbb{E}^\mu\left[\left(\int_{Ns}^{Nt} (\xi_u(0) - \vartheta)\, du\right)^4\right] \leq \begin{cases} CN^2(\log N)^2(t-s)^2, & d=4, \\ CN^2(t-s)^2, & d \geq 5. \end{cases}$$

Obviously, this lemma requires 4th moment computations. We do this by a "bookkeeping of trees" similarly as in [16] for superprocesses. We refer to [11], Lemma 4.1 for details. See also [7], proof of formula (3.36) for a related approach using Laplace transforms.

In the case $d=3$ it turns out that second moments [$\alpha = 2$, $\beta = 1/2$ in (4.1)] suffice. The corresponding estimate is provided in the following lemma.

LEMMA 4.2. *Let $d = 3$ and $\mu \in \{\mathcal{H}(\vartheta), \Lambda_\vartheta\}$. For each $t_0 \geq 0$, there exists a constant $C = C(t_0, \vartheta)$ such that*

$$\mathbb{E}^\mu[(X_t^N - X_s^N)^2] \leq C(t-s)^{3/2} \qquad \forall 0 \leq s \leq t \leq t_0$$

*holds uniformly in $N$.*

Proof. We begin with the independent case. Note that for any initial distribution $\mu$ we have

$$0 \leq \mathbb{E}^\mu[(X_t^N - X_s^N)^2] = \frac{1}{N^{3/2}} \int_{Ns}^{Nt} du \int_{Ns}^{Nt} dv \, \mathrm{Cov}^\mu(\xi_u(0), \xi_v(0)),$$

thus, we see from Lemma 3.3 that $\mathbb{E}^{\mathcal{H}(\vartheta)}[(X_t^N - X_s^N)^2] \leq \mathbb{E}^{\Lambda_\vartheta}[(X_t^N - X_s^N)^2]$ and it is, hence, sufficient to consider the stationary initial distribution $\Lambda_\vartheta$. By stationarity, we can assume without loss of generality that $s = 0$, $t \leq t_0$. Put $\varphi(r) := \mathrm{Cov}^{\Lambda_\vartheta}(\xi_r(0), \xi_0(0))$. We have $0 \leq \varphi(r) \leq C(1 \wedge r^{-1/2})$ by Lemma 3.3 and (2.2). This allows to estimate (and the same is true in the state dependent case)

$$\mathbb{E}^{\Lambda_\vartheta}[(X_t^N - X_0^N)^2]$$
$$= \frac{2}{N^{3/2}} \int_0^{Nt} du \int_u^{Nt} dv \, \varphi(v-u)$$
$$= 2N^{1/2} \int_0^t du \int_u^t dv \, \varphi(N(v-u))$$
$$\leq 2tN^{1/2} \int_0^t dw \, \varphi(Nw)$$
$$\leq 2CN^{1/2} \left\{ t^2 \mathbf{1}(Nt \leq 1) + t\mathbf{1}(Nt > 1)\left[\frac{1}{N} + \int_{1/N}^t ds \frac{1}{\sqrt{Ns}}\right]\right\}$$
$$= 2Ct^{3/2}\left\{(Nt)^{1/2}\mathbf{1}(Nt \leq 1) + (Nt)^{-1/2}\mathbf{1}(Nt > 1)\right.$$
$$\left. + (t/N)^{-1/2}\mathbf{1}(Nt > 1)\left(2(t/N)^{1/2} - \frac{2}{N}\right)\right\}$$
$$\leq 6Ct^{3/2}.$$

For the state dependent case, we observe (again by Lemma 3.3)

$$\mathbb{E}^{\mathcal{H}(\vartheta)}[(X_t^N - X_s^N)^2]$$
$$= \frac{1}{h_d(N)^2} \mathbb{E}^{\mathcal{H}(\vartheta)}\left[\left(\int_{Ns}^{Nt} (\xi_u(0) - \vartheta) \, du\right)^2\right]$$
$$= \frac{1}{h_d(N)^2} \int_{Ns}^{Nt} du \int_{Ns}^{Nt} dv \, (\mathbb{E}^{\mathcal{H}(\vartheta)}[\xi_u(0)\xi_v(0)] - \vartheta^2)$$
$$= \frac{2}{h_d(N)^2} \int_{Ns}^{Nt} du \int_u^{Nt} dv \bigg(\vartheta a_{v-u}(x,y)$$
$$+ \int_0^u a_{v-u+2r}(x,y)\mathbb{E}^{\mathcal{H}(\vartheta)}[\sigma(\xi_{u-r}(0))] \, dr,\bigg)$$

$$\leq \frac{2}{h_d(N)^2} \int_{Ns}^{Nt} du \int_u^{Nt} dv \bigg( \vartheta a_{v-u}(x,y) + c_2 \vartheta \int_0^u a_{v-u+2r}(x,y)\, dr \bigg),$$

which can be estimated as in the independent case. $\square$

**Acknowledgments.** The authors would like to thank Anton Wakolbinger for many stimulating discussions and his constant interest during the preparation of this work. Thanks also to two anonymous referees for their insightful suggestions which helped to streamline the presentation of the paper. Part of this work was completed while the authors enjoyed the hospitality of the Erwin Schrödinger Institute for Mathematical Physics, Vienna.

WEIERSTRASS-INSTITUT
 FÜR ANGEWANDTE ANALYSIS UND STOCHASTIK
MOHRENSTRASSE 39
10117 BERLIN
GERMANY
E-MAIL: birkner@wias-berlin.de

MATHEMATISCHES INSTITUT
UNIVERSITÄT ERLANGEN-NÜRNBERG
BISMARCKSTR. 1 1/2
91054 ERLANGEN
GERMANY
E-MAIL: zaehle@mi.uni-erlangen.de